\documentclass[dvips,ejs,noinfoline]{imsart}

\RequirePackage{amsthm,amsmath,amssymb}

\usepackage{dsfont}
\usepackage{mathbh}
\let\mathds\mathbh

\RequirePackage[colorlinks]{hyperref}

\doi{10.1214/08-EJS288}
\pubyear{2008}
\volume{2}
\issue{0}
\firstpage{1129}
\lastpage{1152}

\startlocaldefs

\numberwithin{equation}{section}

\newtheorem{thm}{Theorem}[section]
\newtheorem{prop}[thm]{Proposition}
\newtheorem{lemma}[thm]{Lemma}
\newtheorem{dfn}{Definition}[section]
\theoremstyle{remark}
\newtheorem{rmk}{Remark}[section]

\begin{document}

\begin{frontmatter}

\title{LASSO, Iterative Feature Selection and the Correlation Selector:
Oracle inequalities and numerical performances}
\runtitle{LASSO, IFS and Correlation Selector}

\begin{aug}
\author{\fnms{Pierre} \snm{Alquier}\ead[label=e1]{alquier@ensae.fr}\ead[label=e2,url]{http://alquier.ensae.net/}\thanksref{t1}}

\thankstext{t1}{I Would like to thank Professors Olivier Catoni, Alexandre Tsybakov, Mohamed Hebiri and Joseph Salmon
as well as the anonymous referees for useful remarks.}

\address{Laboratoire de Probabilit\'es et Mod\`eles Al\'eatoires
(Universit\'e Paris 7)\\
175, rue du Chevaleret\\
75252 Paris Cedex 05, 
France}
\address{CREST, LS \\
3, avenue Pierre Larousse\\ 92240 Malakoff, France\\
\printead{e1}\\ \printead{e2}}

\runauthor{P. Alquier}
\end{aug}


\begin{abstract}
We propose a general family of algorithms for regression estimation
with quadratic loss, on the basis of geometrical considerations.
These algorithms are able to select relevant
functions into a large dictionary. We prove that a lot of methods
that have already been studied for this task (LASSO,
Dantzig selector, Iterative Feature Selection, among others) belong to our family,
and exhibit another particular member of this family that we call
Correlation Selector in this paper.
Using general properties of our family of algorithm we prove
oracle inequalities for IFS, for the LASSO and for the Correlation Selector,
and compare numerical performances of these estimators on a toy example.
\end{abstract}



\received{\smonth{8} \syear{2008}}

\tableofcontents

\end{frontmatter}



\section{Introduction}

\subsection{The regression problem}

In this paper, we study the linear regression problem:
we observe $n$ pairs $(X_{i},Y_{i})$ with $ Y_{i} = f(X_{i}) +
\varepsilon_{i} $ for a noise $\varepsilon =
(\varepsilon_{1},\ldots,\varepsilon_{n})$ to be specified later.

The idea is that the statistician is given (or chooses) a dictionary of functions: $(f_{1},\ldots,f_{m})$,
with possibly $m>n$,
and he wants to build a ``good'' estimation of $f$ of the form $\alpha_{1}f_{1}+\cdots+\alpha_{m}f_{m}$.

Actually, we have to precise two things: what is the distribution of the pairs $(X_{i},Y_{i})$, and
what is the criterion for a ``good'' estimation. We are going to consider two cases.

\subsection{Deterministic and random design}

\subsubsection{Deterministic design case} In this case the values $ X_{1},\ldots, X_{n}$
are deterministic, and the $\varepsilon_{i}$ are i. i. d. according to some distribution
$\mathds{P}$ with $\mathds{E}_{\varepsilon\sim \mathds{P}}(\varepsilon)=0$ and
$\mathds{E}_{\varepsilon\sim \mathds{P}}(\varepsilon^{2})<\infty$.
In this case, the distance between $f$ and $\alpha_{1}f_{1}+\cdots+\alpha_{m}f_{m}$
will be measured in terms of the so-called empirical norm.

\begin{dfn}
For any $\alpha=(\alpha_{1},\ldots,\alpha_{m})\in\mathds{R}^{m}$ and $\alpha'=(\alpha_{1}',\ldots,\alpha_{m}')\in\mathds{R}^{m}$
we put
$$ \left\|\alpha-\alpha'\right\|_{n}^{2} = \frac{1}{n}
\sum_{i=1}^{n}\Biggl[\sum_{j=1}^{m}\alpha_{j}f_{j}(X_{i})-\sum_{j=1}^{m}\alpha_{j}'f_{j}(X_{i})\Biggr]^{2}$$
and
$$ \overline{\alpha}_{n} \in \arg\min_{\alpha\in\mathds{R}^{m}}
         \frac{1}{n} \sum_{i=1}^{n}\Bigg[f(X_{i})-\sum_{j=1}^{m}\alpha_{j}f_{j}(X_{i}) \Bigg]^{2}.$$
\end{dfn}

\subsubsection{Random design case} In this case, we assume that the pairs $(X_{i},Y_{i})$ are i. i. d. according
to some distribution $\mathds{P}$,
that the marginal distribution of every $X_{i}$ is $\mathds{P}_{X}$,
and that we still have $\mathds{E}_{(X,Y)\sim \mathds{P}}(\varepsilon)=0$ and
$\mathds{E}_{(X,Y)\sim \mathds{P}}(\varepsilon^{2})<\infty$. The distance will be measured by
the $\mathfrak{L}^{2}$ distance with respect to $\mathds{P}_{X}$.

\begin{dfn}
For any $\alpha,\alpha'\in\mathds{R}^{m}$
we put
$$ \left\|\alpha-\alpha'\right\|_{X}^{2} = \mathds{E}_{X\sim \mathds{P}_{X}}
\left\{\Bigg[\sum_{j=1}^{m}\alpha_{j}f_{j}(X)-\sum_{j=1}^{m}\alpha_{j}'f_{j}(X) \Bigg]^{2}\right\}$$
and
$$ \overline{\alpha}_{X} \in \arg\min_{\alpha\in\mathds{R}^{m}}
          \mathds{E}_{X\sim \mathds{P}_{X}}
              \left\{\Bigg[f(X)-\sum_{j=1}^{m}\alpha_{j}f_{j}(X) \Bigg]^{2}\right\}.$$
\end{dfn}

Moreover, we make the following restrictive hypothesis: the statistician\break knows~$\mathds{P}_{X}$.

\subsection{General notations}

Now, we assume that we are in one of the two cases defined previously. However, as the results we want
to state are the same in both settings, we introduce the following notation.

\begin{dfn}
We introduce the general norm
$$ \|\alpha-\alpha'\|_{GN} $$
that is simply $\left\|\alpha-\alpha'\right\|_{n}$ if we are in the deterministic design case
and and $ \left\|\alpha-\alpha'\right\|_{X}$ if we are in the random design case.
Moreover, we will let $\overline{\alpha}$ denote $\overline{\alpha}_{n}$ or $\overline{\alpha}_{X}$
according to the case.

In any case, we let $P$ denote the distribution of the sample $(X_{i},Y_{i})_{i=1,\ldots,n}$.
\end{dfn}

In order so simplify the notations, we assume that the functions $f_{j}$ of the dictionary are normalized,
in the sense that $ \frac{1}{n}\sum_{i=1}^{n} f_{j}^{2}(X_{i}) = 1 $
if we are in the deterministic design case and that $ \mathds{E}_{X\sim \mathds{P}_{X}}\left[f_{j}(X)\right]^{2}=1$
if we are in the random design case. Note that this could be simply written in terms of the general norm:
if we put $e_{1}=(1,0,\ldots,0)$,\ldots, $e_{m}=(0,\ldots,0,1)$ the canonical basis of $\mathds{R}^{m}$, we just have
to assume that for any $j\in\{1,\ldots,m\}$, $\|e_{j}\|_{GN}=1$.

Finally, let us mention that $\left<.,.\right>_{GN}$ will denote the scalar product associated to the
norm $\|.\|_{GN}$ while we will use the notation $\|.\|$ for the euclidian norm in $\mathds{R}^{m}$ and
$\left<.,.\right>$ for the associated scalar product.

\subsection{Previous works and organization of the paper}

The aim of this paper is to propose a method to estimate the real
regression function (say $f$) on the basis of the dictionary $(f_{1},\ldots,f_{m})$,
that have good performances even if $m>n$.

Recently, a lot of algorithms have been proposed for that purpose,
let's cite among others the bridge regression by Frank and Friedman
\cite{bridge}, and a particular case of bridge regression called LASSO by Tibshirani \cite{Lasso1},
some variants or generalization like LARS by Efron, Hastie,
Johnstone and Tibshirani \cite{LARS}, the Dantzig selector by Candes
and Tao \cite{Dantzig} and the Group LASSO by Bakin \cite{Bakin}, Yuan and Lin
\cite{GLasso} and Chesneau and Hebiri \cite{GVLasso} or iterative algorithms
like Iterative Feature Selection in our paper
\cite{Alquier2007} or greedy algorithms in Barron, Cohen, Dahmen and DeVore \cite{greedy}.
This paper proposes a general method that contains LASSO, Dantzig selector and
Iterative Feature Selection as a particular case.

Note that in the case where $m/n$ is small, we can use the ordinary
least square estimate. The risk of this estimator is roughly in
$m/n$. But when $m/n >1$, this estimator isn't even properly
defined. The idea of all the mentioned works is the following: if
there is a ``small'' vector space $F\subset\mathds{R}^{m}$ such that
$\overline{\alpha}\in F$, one could build a constrained estimator
with a risk in $\dim(F)/n$. But can we obtain such a result if $F$
is unknown? For example, a lot of papers study the {\it sparsity} of
$\overline{\alpha}$, this means that $F$ is the span of a few
$e_{j}$, or, in other words, that $\overline{\alpha}$ have only a
small number (say $p$) of non-zero coordinates: an estimator that
selects automatically $p$ relevant coordinates and achieving a
risk close to $p/n$ is said to satisfy a ``sparsity oracle
inequality''. A paper, by Bickel, Ritov and Tsybakov \cite{Lasso3}
gives sparsity oracle inequalities for the LASSO and the Dantzig
selector in the case of the deterministic design. Another paper by
Bunea, Tsybakov and Wegkamp \cite{Lasso2} gives sparsity oracle
inequalities for the LASSO. This paper is written in a more general
context than ours: random design with {\it unknown} distribution (in
the case of a random design, remember that our method require the
knowledge of the distribution of the design). However, the main
results require the assumption $\|f_{j}\|_{\infty}\leq L$ for some
given $L$, what is not necessary in our paper, and prevents the
use of popular basis of functions like wavelets. This is due to the
use of Hoeffding's inequality in the technical parts of the
paper.\looseness=1

Our paper uses a geometric point of view. This allows to build a
general method of estimation and to obtain simple sparsity oracle
inequalities for the obtained estimator, in both deterministic
design case and random design with known distribution. It uses a
(Bernstein's type) deviation inequality proved in a previous work
\cite{Alquier2007} that is sharper than Hoeffding's inequality, and
so gets rid of the assumption of a (uniform) bound over the functions
of the dictionary. Another improvement is that our method is valid
for some types of data-dependant of dictionaries of functions, for
example the case where $m=n$ and
$$ \{f_{1}(.),\ldots,f_{m}(.)\}=\{K(X_{1},.),\ldots,K(X_{n},.)\} $$
where $K$ is a function $\mathcal{X}^{2}\rightarrow\mathds{R}$, performing kernel estimation.

In Section \ref{general}, we give the general form for our algorithm
under a particular assumption, Assumption {\bf (CRA)}, that says we are able to build some
confidence region for the best value of $\alpha$ in some subspace of
$\mathds{R}^{m}$.

In Section \ref{oracle}, we show why Iterative Feature Selection
(IFS), LASSO, Dantzig Selector among others are particular cases of
our algorithm. We exhibit another particular case of interest
(called the Correlation Selector in this paper). Moreover, we prove
some oracle inequalities for the obtained estimators: roughly,
LASSO, Dantzig Selector and IFS performs well when the vector
$\overline{\alpha}$ is sparse (which means that a lot of its
coordinates,
$\left<\overline{\alpha},e_{j}\right>=\overline{\alpha}_{j}$ are
equal to zero) or approximately sparse (a lot of coordinates are
nearly equal to zero), while the Correlation Selector performs well
when a lot of $\left<\overline{\alpha},e_{j}\right>_{GN}$ are almost
equal to zero (in the deterministic design case,
$\left<\overline{\alpha},e_{j}\right>_{GN}=\mathds{E}(\frac{1}{n}\sum_{i=1}^{n}f_{j}(X_{i})Y_{i})$
while in the random design case,
$\left<\overline{\alpha},e_{j}\right>_{GN}=\mathds{E}(f_{j}(X)Y)$,
so in any case, this quantity is a measure of the correlation
between the variable $Y$ and the $j$-th function in the dictionary).
 So, intuitively, the Correlation Selector gives good results when most
 of the functions in the dictionary have weak
correlation with $Y$, but we expect that altogether these functions
can bring a good prediction for $Y$.

In order to prove oracle inequalities, some types of orthogonality
(or approximate orthogonality, in some sense) are required on the
dictionary of functions. Our results are the following: under
orthogonality on the dictionary of functions, and {\it using only
general properties of our family of estimators}, we have a sparse
oracle inequality. Under an approximate orthogonality condition
taken from Bickel, Ritov and Tsybakov \cite{Lasso3}, the result can
be extended for the LASSO and the Dantzig selector (with a proof
taken from \cite{Lasso3}). Some remarks by Huang, Cheang and Barron
\cite{Barron2} show that these results can be extended to IFS with a
slight modification of the estimator. Finally, the central result for
the Correlation Selector does not require any hypothesis on the
dictionary of functions but concerns a measure of the risk that is
not natural, we obtain a result on the risk measured by $\|.\|_{GN}$
under an assumption very close to the one in \cite{Lasso3} - here
again, the proof uses only general properties of our family of
estimators.

Section \ref{simulations} is dedicated to simulations:
we compare ordinary least square (OLS), LASSO, Iterative Feature
Selection and the Correlation Selector on a toy example. Simulations shows that
both particular cases of our family of estimators (LASSO and Iterative Feature Selection)
generally outperforms the OLS estimate. Moreover, LASSO performs
generally better than Iterative Feature Selection, however, this is
not always true: this fact leads to the conclusion that a data-driven
choice of a particular algorithm in our general family could lead
to optimal results.

After a conclusion (Section \ref{secconc}), Section \ref{proofs} is dedicated to some proofs.

\section{General projection algorithms}

\label{general}

\subsection{Additional notations and hypothesis}

\begin{dfn}
Let $\mathcal{C}$ be a closed, convex subset of $\mathds{R}^{d}$. We
let $\Pi^{GN}_{\mathcal{C}}(.)$ denote the orthogonal projection on
$\mathcal{C}$ with respect to the norm $\|.\|_{GN}$:
$$ \Pi^{GN}_{\mathcal{C}}(\alpha) = \arg\min_{\beta\in\mathcal{C}} \|\alpha-\beta\|_{GN}. $$
For a generic distance $\delta$, we will use the notation $\Pi^{\delta}_{\mathcal{C}}(.)$
for the orthogonal projection on $\mathcal{C}$ with respect to $\delta$.
\end{dfn}

We put, for every $j\in\{1,\ldots,m\}$:
$$ \mathcal{M}_{j} = \bigl\{\alpha\in\mathds{R}^{m},\quad \ell\neq j\Rightarrow \alpha_{\ell}=0 \bigr\}
                   = \left\{\alpha e_{j},\alpha\in\mathds{R}\right\}. $$

\begin{dfn}
We put, for every $j\in\{1,\ldots,m\}$:
$$ \overline{\alpha}^{j} = \arg\min_{\alpha\in\mathcal{M}_{j}} \|\overline{\alpha}-\alpha e_{j} \|_{GN}
= \Pi^{GN}_{\mathcal{M}_{j}}(\overline{\alpha}). $$
Moreover let us put:
$$ \tilde{\alpha}_{j}= \frac{1}{n}\sum_{i=1}^{n} f_{j}(X_{i}) Y_{i} \text{ and }
         \hat{\alpha}^{j} =  \tilde{\alpha}_{j}  e_{j}. $$
\end{dfn}

\begin{rmk}
In the deterministic design case ($\|.\|_{GN}=\|.\|_{n}$) we have
$$ \overline{\alpha}^{j} = \Bigg[\frac{1}{n}\sum_{i=1}^{n} f_{j}(X_{i}) f(X_{i})\Bigg] e_{j} $$
and in the random design case we have
$$ \overline{\alpha}^{j} = \mathds{E}_{X\sim \mathds{P}_{X}} \left[f_{j}(X) f(X)\right] e_{j} $$
so in any case, $\hat{\alpha}^{j}$ is an estimator of $\overline{\alpha}^{j}$.
\end{rmk}

\noindent {\bf Hypothesis (CRA)}
We say that the confidence region assumption (CRA) is
satisfied if for $\varepsilon\in[0,1]$ we have a bound
$r(j,\varepsilon)\in\mathds{R}$ such that
$$ P\Bigl[\forall j\in\{1,\ldots,m\},\quad \left\|\overline{\alpha}^{j}-\hat{\alpha}^{j}\right\|_{X}^{2}
         \leq r(j,\varepsilon) \Bigr]
 \geq 1-\varepsilon. $$

In our previous work \cite{Alquier2007} we examined different
hypothesis on the probability $P$ such that this hypothesis is
satisfied. For example, using inequalities by Catoni \cite{Catoni} and
Panchenko \cite{Panchenko} we
proved the following results.

\begin{lemma}
Let us assume that $\|f\|_{\infty} \leq L$ for some known $L$.
Let us assume that $\mathds{E}_{\mathds{P}}(\varepsilon^{2})\leq \sigma^{2}$ for some known $\sigma^{2}<\infty$.
Then Assumption {\bf (CRA)} is satisfied, with
$$ r(j,\varepsilon) = \frac{4\left(1+\log\frac{2m}{\varepsilon}\right)}{n}
     \Bigg[\frac{1}{n}\sum_{i=1}^{n}f_{j}^{2}(X_{i})Y_{i}^{2} + L^{2}+\sigma^{2} \Bigg] .$$
\end{lemma}

\begin{rmk}
It is also shown in \cite{Alquier2007} that we are allowed to take
$$ \{f_{1}(.),\ldots,f_{m}(.)\} =  \{K(X_{1},.),\ldots,K(X_{n},.)\} $$
for some function $K: \mathcal{X}^{2}\rightarrow\mathds{R}$
(this allows for $f(x)$ a kernel estimator of the form $\sum_{i=1}^{n}\alpha_{i}K(X_{i},x)$),
even in the random design case, but we have to take
$$ r(j,\varepsilon) = \frac{4\left(1+\log\frac{4m}{\varepsilon}\right)}{n}
     \Bigg[\frac{1}{n}\sum_{i=1}^{n}f_{j}^{2}(X_{i})Y_{i}^{2} + L^{2}+\sigma^{2} \Bigg] $$
in this case.
\end{rmk}

\begin{lemma}
Let us assume that there is a $K>0$ such that $\mathds{P}(|Y|\leq K)=1$.
Then Assumption {\bf (CRA)} is satisfied with
$$ r(j,\varepsilon) = \frac{8K^{2}\left(1+\log\frac{2m}{\varepsilon}\right)}{n}.$$
\end{lemma}

\begin{dfn}
When {\bf (CRA)} is satisfied,
we define, for any $\varepsilon>0$ and $j\in\{1,\ldots,m\}$, the random set
$$ \mathcal{CR}(j,\varepsilon) = \Bigl\{\alpha\in\mathds{R}^{m},
 \quad \left\|\Pi^{GN}_{\mathcal{M}_{j}}(\alpha)-\hat{\alpha}^{j}\right\|^{2}_{GN} \leq r(j,\varepsilon)\Bigr\}. $$
\end{dfn}

This can easily be interpreted: Assumption {\bf (CRA)} says that
there is a confidence region for $\overline{\alpha}^{j}$ in the
small model $\mathcal{M}_{j}$; $ \mathcal{CR}(j,\varepsilon)$ is the
set of all vectors falling in this confidence region when they are
orthogonally  projected on $\mathcal{M}_{j}$.

We remark that the hypothesis implies that
$$ P\Bigl[\forall j\in\{1,\ldots,M\},\quad \overline{\alpha} \in \mathcal{CR}(j,\varepsilon) \Bigr] \geq 1-\varepsilon. $$

\subsection{General description of the algorithm}

We propose the following iterative algorithm.
Let us choose a confidence level $\varepsilon>0$ and
a distance on $\mathcal{X}$, say $\delta(.,.)$.

\begin{itemize}
\item Step 0. Choose $\hat{\alpha}(0)=(0,\ldots,0)\in\mathds{R}^{m}$.
Choose $\varepsilon\in[0,1]$.
\item General Step ($k$). Choose $N(k)\leq M$ and indices
$(j_{1}^{(k)},\ldots,j_{N}^{(k)})\in\break \{1,\ldots,M\}^{N(k)}$ and put:
$$ \hat{\alpha}(k) \in \arg\min_{\alpha\in\bigcap_{\ell=1}^{N(k)}\mathcal{CR}(j^{(k)}_{\ell},\varepsilon)}
\delta\bigl(\alpha,\hat{\alpha}(k-1)\bigr).$$
\end{itemize}

This algorithm is motivated by the following result.

\begin{thm}
\label{THM} When the CRA assumption is satisfied we have:
\begin{equation}
\label{MOTIV} P\Big[\forall k\in\mathds{N},\quad
\delta\bigl(\hat{\alpha}(k),\overline{\alpha}\bigr) \leq
 \delta\bigl(\hat{\alpha}(k-1),\overline{\alpha}\bigr) \leq \cdots
\leq \delta\bigl(\hat{\alpha}(0),\overline{\alpha}\bigr)
\Big]\geq 1-\varepsilon .
\end{equation} So, our algorithm builds a
sequence of $\hat{\alpha}(k)$ that gets closer to
$\overline{\alpha}$ (according to $\delta$) at every step. Moreover,
if $\delta(x,x')=\left\|x-x'\right\|_{GN}$ then
$$ \hat{\alpha}(k) = \Pi^{GN}_{\bigcap_{\ell=1}^{N(k)}\mathcal{CR}(j^{(k)}_{\ell},\varepsilon) }
                 \bigl(\hat{\alpha}(k-1)\bigr) $$
and we have the following:
$$ P\Biggl[\forall k\in\mathds{N}, \quad \left\|\hat{\alpha}(k)-\overline{\alpha}\right\|_{GN}^{2}
\leq
\left\|\hat{\alpha}(0)-\overline{\alpha}\right\|_{GN}^{2}
 - \sum_{j=1}^{k}\left\|\hat{\alpha}(j)-\hat{\alpha}(j-1)\right\|_{GN}^{2} \Bigg]\,{\geq}\, 1-\varepsilon .$$
\end{thm}

\begin{proof}
Let us assume that
$$ \forall j\in\{1,\ldots,M\},\quad \left\|\overline{\alpha}^{j}-\hat{\alpha}^{j}\right\|_{GN}^{2}
\leq r(S_{j},\varepsilon). $$
This is true with probability at least $1-\varepsilon$ according to assumption {\it (CRA)}.
In this case we have seen that
$$\overline{\alpha} \in \bigcap_{\ell=1}^{N(k)}\mathcal{CR}(j^{(k)}_{\ell},\varepsilon) $$
that is a closed convex region,
and so, by definition, $\delta\left(\hat{\alpha}(k),\overline{\alpha}\right) \leq
 \delta\left(\hat{\alpha}(k-1),\overline{\alpha}\right)$ for any $k\in\mathds{N}$.
If $\delta$ is the distance associated with the norm $\|.\|_{GN}$,
let us choose $k\in\mathds{N}$,
\begin{multline*}
\left\|\hat{\alpha}(k)-\overline{\alpha}\right\|_{GN}^{2}
=  \left\|\Pi^{GN}_{\bigcap_{\ell=1}^{N(k)}\mathcal{CR}(j^{(k)}_{\ell},\varepsilon)}
\bigl(\hat{\alpha}(k-1)\bigr)-\overline{\alpha}\right\|_{GN}^{2}
\\
\leq \left\|\hat{\alpha}(k-1)-\overline{\alpha}\right\|_{GN}^{2}
  - \left\|\Pi^{GN}_{\bigcap_{\ell=1}^{N(k)}\mathcal{CR}(j^{(k)}_{\ell},\varepsilon)}
            \bigl(\hat{\alpha}(k-1)\bigr)-\hat{\alpha}(k-1)\right\|_{GN}^{2}
\\
= \left\|\hat{\alpha}(k-1)-\overline{\alpha}\right\|_{GN}^{2} - \left\|\hat{\alpha}(k)-\hat{\alpha}(k-1)\right\|_{GN}^{2}.
\end{multline*}
A recurrence ends the proof.
\end{proof}

\begin{rmk}
We choose our estimator $\hat{\alpha}=\hat{\alpha}(k)$ for some
step $k\in\mathds{N}$; the choice of the stopping step $k$ will
depend on the particular choices of the projections and is detailed
in what follows. But remark that \textbf{there is no bias-variance
balance involved in the choice of} $k$ as Theorem \ref{THM} shows
that overfitting is not possible for large values of $k$.
\end{rmk}

\section{Particular cases and oracle inequalities}

\label{oracle}

We study some particular cases depending on the choice of the distance $\delta(.,.)$ and
on the sets we project on.

Roughly, LASSO and Iterative Feature Selection (at least as introduced in \cite{Alquier2007})
correspond to the choice $\delta(\alpha,\alpha')=\|\alpha-\alpha'\|_{GN}$, and are studied first.

Dantzig selector corresponds to the choice $\delta(\alpha,\alpha')=\|\alpha-\alpha'\|_{1}$ the $\ell_{1}$ distance,
it is studied in a second time.

Finally, the new Correlation Selector corresponds to another choice for $\delta$.

\subsection{The LASSO}

Here, we use only one step where we project $0$ onto the
intersection of all the confidence regions and so we obtain:
$$  \hat{\alpha}^{L} = \hat{\alpha}(1)
= \Pi^{GN}_{\bigcap_{\ell=1}^{m}\mathcal{CR}(\ell,\varepsilon) }\left(0\right) . $$

The optimization program to obtain $\hat{\alpha}^{L}$ is given by:
$$
\left\{
\begin{array}{l}
\arg\min_{\alpha=(\alpha_{1},\ldots,\alpha_{m})\in\mathds{R}^{m}}
\|\alpha\|_{GN}^{2}
\\
\\[-6pt]
\text{s. t.} \quad
\alpha\in\bigcap_{\ell=1}^{m}\mathcal{CR}(\ell,\varepsilon)
\end{array}
\right.
$$
and so:
\begin{equation}
\label{ProgLASSO}
\left\{
\begin{array}{l}
\arg\min_{\alpha\in\mathds{R}^{m}} \|\alpha\|_{GN}^{2}
\\
\\[-6pt]
\text{s. t.} \quad \forall j\in\{1,\ldots,m\},\quad
\left|\left<\alpha,e_{j}\right>_{GN}-\tilde{\alpha}_{j}\right| \leq
\sqrt{r(j,\varepsilon)}
\end{array}
\right.
\end{equation}

\begin{prop}
\label{PropLASSO} Every solution of the program
\begin{equation}
\label{lassogeneral} \arg\min_{\alpha\in\mathds{R}^{m}}
\left\{\|\alpha\|_{GN}^{2}-2\sum_{j=1}^{m}\alpha_{j}\tilde{\alpha}_{j}
+ 2\sum_{j=1}^{m} \sqrt{r(j,\varepsilon)}
\left|\alpha_{j}\right|\right\}
\end{equation}
satisfies Program \ref{ProgLASSO}.
Moreover, all the solutions $\alpha$ of Program \ref{ProgLASSO} have the same risk
   value $\|\alpha-\overline{\alpha}\|_{GN}^{2}$.
Finally, in the \textbf{deterministic
design case}, Program \ref{lassogeneral} is equivalent to:
\begin{equation}
\label{las2}
\arg\min_{\alpha\in\mathds{R}^{m}} \Bigg\{
\frac{1}{n}\sum_{i=1}^{n}\biggl[Y_{i}-\sum_{j=1}^{m}\alpha_{j} f_{j}(X_{i})\biggr]^{2}
+ 2\sum_{j=1}^{m} \sqrt{r(j,\varepsilon)}\left|\alpha_{j}\right|\Bigg\}.
\end{equation}
\end{prop}

The proof is given ay the end of the paper (in Subsection \ref{ProofLASSO} page \pageref{ProofLASSO}).

Note that, if $r(j,\varepsilon)$ does not depend on $j$, Program \ref{las2} is
exactly one of the formulations of the LASSO estimator studied first by
Tibshirani \cite{Lasso1}. In the particular \textbf{deterministic
design case}, this dual representation was already known and introduced by Osborne, Presnell and
Turlach \cite{lassodual}.

However, in the cases where $r(j,\varepsilon)$ is not constant, the
difference with the LASSO algorithm is the following: the harder the coordinates are to be
estimated, the more penalized they are.

Moreover, note that the program \ref{lassogeneral} gives a
different from of the usual LASSO program for the cases where we do not
use the empirical norm.

\subsection{Iterative Feature Selection (IFS)}

\label{IFSdes}

As in the LASSO case we use the distance
$\delta(\alpha,\beta)=\|\alpha-\beta\|_{GN}$.

The only difference is that instead of taking the intersection of
every confidence region, we project on each of them iteratively. So
the algorithm is the following:
$$  \hat{\alpha}(0) = (0,\ldots,0) $$
and at each step $k$ we choose a $j(k)\in\{1,\ldots,m\}$ and
$$  \hat{\alpha}(k) = \Pi^{GN}_{\mathcal{CR}(j(k),\varepsilon) }\bigl(\hat{\alpha}(k-1)\bigr). $$
We choose a stoping step $\hat{k}$ and put
$$ \hat{\alpha}^{IFS}= \hat{\alpha}(\hat{k}) .$$

This is exactly the Iterative Feature Selection algorithm that was
introduced in Alquier \cite{Alquier2007}, with the choice of $j(k)$:
$$ j(k) = \arg\max_{j} \left\| \hat{\alpha}(k-1)
- \Pi^{X}_{\mathcal{CR}(j,\varepsilon) }\bigl(\hat{\alpha}(k-1)\bigr) \right\|_{GN} ,$$
and the suggestion to take as a stopping step
$$ \hat{k} = \inf\left\{k\in\mathds{N}^{*},\quad \left\|  \hat{\alpha}(k)
   -\hat{\alpha}(k-1) \right\|_{GN} \leq \kappa\right\} $$
for some small $\kappa>0$.

\begin{rmk}
In Alquier \cite{Alquier2007}, it is proved that:
\begin{equation}
\label{IFS1}
  \hat{\alpha}(k) = \hat{\alpha}(k-1) + sgn\left(\beta_{k}\right)
            \left(\left|\beta_{k}\right|-\sqrt{r(j(k),\varepsilon)}\right)_{+} e_{j(k)}
\end{equation}
where
$$ \beta_{k} = \frac{1}{n}\sum_{i=1}^{n} f_{j}(X_{i})\Bigg[Y_{i}
                   - \sum_{\ell=1}^{m}\hat{\alpha}(k)_{\ell}f_{\ell}(X_{i})\Bigg] .$$
                   So this algorithm looks quite similar to a greedy algorithm,
as it is described by Barron, Cohen, Dahmen and DeVore
\cite{greedy}. Actually, it would be a greedy algorithm if we
replace $r(j,\varepsilon)$ by $0$ (such a choice is however not
possible here): it is a soft-thresholded version of a greedy
algorithm. Such greedy algorithms were studied in a recent
paper by Huang, Cheang and Barron \cite{Barron2} under the name
``penalized greedy algorithm'', in the case $\|.\|_{GN}=\|.\|_{n}$.

Note that in Iterative
Feature Selection, every selected feature actually improves the
estimator:  $\|\hat{\alpha}(k)-\overline{\alpha}\|_{GN}^{2}\leq \|\hat{\alpha}(k-1)-\overline{\alpha}\|_{GN}^{2}$
(Equation \ref{MOTIV}).
\vspace*{-6pt}
\end{rmk}

\subsection{The Dantzig selector}

The Dantzig selector is based on a change of distance $\delta$. We choose
$$\delta(\alpha,\alpha')=\|\alpha-\alpha'\|_{1}=\sum_{j=1}^{m}|\alpha_{j}-\alpha_{j}'|.$$

As is the LASSO case, we make only one projection onto the
intersection of every confidence region:
$$  \hat{\alpha}^{DS} \in
\arg\min_{\alpha\in\bigcap_{\ell=1}^{m}\mathcal{CR}(j,\varepsilon)}
\|\alpha\|_{1}
$$
and so $\hat{\alpha}^{DS}$ is the solution of the program:
$$
\left\{
\begin{array}{l}
\displaystyle\arg\min_{\alpha=(\alpha_{1},\ldots,\alpha_{m})\in\mathds{R}^{m}} \sum_{j=1}^{m}\left|\alpha_{j}\right|
\\
\\[-6pt]
\text{s. t.} \quad \forall j\in\{1,\ldots,m\},\quad
\left|\left<\alpha,e_{j}\right>_{GN}-\tilde{\alpha}_{j}\right| \leq
\sqrt{r(j,\varepsilon)}.
\end{array}
\right.
$$
In the case where $r(j,\varepsilon)$ does not depend on $j$, and where $\|.\|_{GN}=\|.\|_{n}$, this
program is exactly the one proposed by Candes and Tao \cite{Dantzig} when they introduced the Dantzig selector.

\subsection{Oracle Inequalities for the LASSO, the Dantzig Selector and IFS}

\begin{dfn}
For any $S\subset\{1,\ldots,m\}$ let us put
$$ \mathcal{M}_{S} = \left\{\alpha\in\mathds{R}^{m},\quad j\notin S\Rightarrow \alpha_{j} = 0\right\} $$
and
$$ \overline{\alpha}_{S} = \arg\min_{\alpha\in\mathcal{M}_{S}} \|\alpha-\overline{\alpha}\|_{GN}. $$
\end{dfn}

Every $ \mathcal{M}_{S} $ is a submodel of $\mathds{R}^{m}$ of dimension $|S|$ and $
\overline{\alpha}_{S}$ is the best approximation of
$\overline{\alpha}$ in this submodel.

\begin{thm}
\label{thm3} Let us assume that assumption {\bf (CRA)} is satisfied. Let
us assume that the functions $f_{1},\ldots,f_{m}$ are orthogonal with
respect to $\left<.,.\right>_{GN}$. In this case the order of the
projections in Iterative Feature Selection does not affect the
obtained estimator, so we can set
$$\hat{\alpha}^{IFS}=\Pi^{GN}_{\mathcal{CR}(m,\varepsilon) } \ldots \Pi^{GN}_{\mathcal{CR}(1,\varepsilon) } 0  .$$
Then
$$ \hat{\alpha}^{IFS} = \hat{\alpha}^{L} = \hat{\alpha}^{DS}
 = \sum_{j=1}^{m} sgn \left(\tilde{\alpha}_{j}\right) \left(\left|\tilde{\alpha}_{j}\right|-\sqrt{r(j,\varepsilon)}\right)_{+} e_{j} $$
is a soft-thresholded estimator, and
$$
P\left\{
\left\|\hat{\alpha}^{L}-\overline{\alpha}\right\|^{2}_{GN}
\leq \inf_{S\subset\{1,\ldots,m\}} \Bigg[
\left\|\overline{\alpha}_{S}-\overline{\alpha}\right\|^{2}_{GN}
+ 4 \sum_{j\in S} r(j,\varepsilon) \Bigg]\right\}
\geq 1-\varepsilon.
$$
\end{thm}

For the proof, see Subsection \ref{proofthm3} page \pageref{proofthm3}.

\begin{rmk}
We call ``general regularity assumption with order $\beta>0$
and constant $C>0$'':
$$ \forall j\in\{1,\ldots,m\},\quad \inf_{
\tiny{
\begin{array}{c}
S\subset\{1,\ldots,m\}
\\
|S|\leq j
\end{array}
} } \left\|\overline{\alpha}_{S} -\overline{\alpha} \right\|_{GN}
\leq C j^{-\beta} .$$
This is the kind of regularity satisfied by functions in weak Besov spaces, see
Cohen \cite{Cohen} and the references therein, with $f_{j}$ being wavelets.
If the general
regularity assumption is satisfied
with regularity $\beta>0$ and constant $C>0$ and if there is a $k>0$
such that
$$ r(j,\varepsilon) \leq \frac{k \log \frac{m}{\varepsilon}}{n}, $$ then we have:
$$
P \left\{ \left\|\hat{\alpha}^{L}-\overline{\alpha}\right\|^{2}_{GN} \leq
\left(2\beta+1\right)C^{\frac{1}{2\beta+1}}\left(\frac{2
k \log \frac{m}{\varepsilon}}{\beta
n}\right)^{\frac{2\beta}{2\beta+1}}+
\left(\frac{4 k \log
\frac{m}{\varepsilon}}{n}\right)\right\}
\geq 1-\varepsilon.
$$
\end{rmk}

Now, note that the orthogonality assumption is very restrictive.
Usual results about LASSO or Dantzig Selector usually involve only
approximate orthogonality, see for example Candès and Tao
\cite{Dantzig}, Bunea \cite{Bunea}, Bickel, Ritov and Tsybakov \cite{Lasso3} and Bunea,
Tsybakov and Wegkamp \cite{Lasso2}, and sparsity (the fact that a
lot of the coordinates of $\overline{\alpha}$ are null), as for example the
following result, which is a small variant of a result in \cite{Lasso2}, that is reminded here
in order to provide comparison with the results coming later in
the paper.

\begin{thm}[{\rm Variant of Bunea, Tsybakov and Wegkamp \cite{Lasso2}}]
\label{thTSYBA}
Let us assume that we are in the \textbf{deterministic design case},
that assumption {\bf (CRA)} is satisfied, and that
$r(j,\varepsilon)=r(\varepsilon)$ does not depend on $j$ (this is always possible by taking
$r(\varepsilon)
=\sup_{j\in\{1,\ldots,m\}}r(j,\varepsilon)$). Moreover, we assume that
there is a constant $D$ such that, for any $\alpha\in\mathcal{F}_{m}$,
$$ \|\alpha\|_{GN} \geq D\|\alpha\| $$
where
$$ \mathcal{F}_{m} = \Bigg\{\alpha\in\mathds{R}^{m},\quad \sum_{j:\overline{\alpha}_{j} = 0}
                            |\alpha_{j}| \leq 3 \sum_{j:\overline{\alpha}_{j} \neq 0}
                            |\alpha_{j}| \Bigg\} .$$
Then
\begin{equation*}
P\Biggl\{
 \left\|\hat{\alpha}^{L}-\overline{\alpha}\right\|^{2}_{GN}
  \leq \frac{16}{D^{2}} \Bigl|\left\{j:\overline{\alpha}_{j}\neq 0\right\}\Bigr|
            r(\varepsilon) \Biggr\} \geq 1-\varepsilon.
\end{equation*}
\end{thm}

The only difference with the original result in \cite{Lasso2} is
that $r(\varepsilon)$ is given in a general form here, so we are
allowed to use different values for $r(\varepsilon)$ depending on
the context, see the discussion of Hypothesis {\bf (CRA)} in the
beginning of the paper. Similar results are available for the
Dantzig selector, see Candes and Tao \cite{Dantzig}, and Bickel,
Ritov and Tsybakov \cite{Lasso3}.

\begin{rmk}
We can wonder how IFS performs when the dictionary is not orthogonal.
Actually, the study of penalized greedy algorithm in Huang, Cheang and Barron \cite{Barron2}
leads to the following conclusion in the deterministic design case: there are cases where
IFS can be really worse than LASSO. However, the authors proposes a modification of the algorithm,
called ``relaxed penalized greedy algorithm''; if we apply this modification here we obtain
\begin{equation*}
  \hat{\alpha}(k) = \gamma_{k} \hat{\alpha}(k-1) + sgn\left(\beta_{k}\right)
            \left(\big|\tilde{\beta}_{k}\big|-\sqrt{r(j(k),\varepsilon)}\right)_{+} e_{j(k)}
\end{equation*}
instead of equation \ref{IFS1}, where
$$ \tilde{\beta}_{k} = \frac{1}{n}\sum_{i=1}^{n} f_{j}(X_{i})\Bigg[Y_{i}
                   - \gamma_{k}\sum_{\ell=1}^{m}\hat{\alpha}(k)_{\ell}f_{\ell}(X_{i})\Bigg],$$
and at each step we have to minimize the empirical least square
error with respect to $\gamma_{k}\in[0,1]$. Such a modification
ensures that the estimators given by the $k$-th step of the
algorithm become equivalent to the LASSO when $k$ grows, for more
details see \cite{Barron2} (note that the interpretation in terms of
confidence regions and the property
$\|\hat{\alpha}(k)-\overline{\alpha}\|_{GN}^{2}\leq
\|\hat{\alpha}(k-1)-\overline{\alpha}\|_{GN}^{2}$ are lost with this
modification).
\end{rmk}

\subsection{A new estimator: the Correlation Selector}

The idea of the Correlation Selector is to use
$$ \left\|\alpha \right\|_{CS}
        =  \sum_{j=1}^{m}\left<e_{j},\alpha\right>_{GN}^{2} .$$

We make only one projection onto the intersection of every
confidence region:
$$  \hat{\alpha}^{CS} \in
\arg\min_{\alpha\in\bigcap_{\ell=1}^{m}\mathcal{CR}(j,\varepsilon)}
\|\alpha\|_{CS}
$$
and so $\hat{\alpha}^{CS}$ is a solution of the program:
$$
\left\{
\begin{array}{l}
\displaystyle\arg\min_{\alpha=(\alpha_{1},\ldots,\alpha_{m})\in\mathds{R}^{m}}
\sum_{j=1}^{m}\left<e_{j},\alpha\right>_{GN}^{2}
\\
\\[-9pt]
\text{s. t.} \quad \forall j\in\{1,\ldots,m\},\quad
\left|\left<\alpha,e_{j}\right>_{GN}-\tilde{\alpha}_{j}\right| \leq
\sqrt{r(j,\varepsilon)}.
\end{array}
\right.
$$
This program can be solved for every $u_{j}=\left<e_{j},\alpha\right>_{GN}$ individually: each of them is solution
of
$$
\left\{
\begin{array}{l}
\arg\min_{u} |u|^{2}
\\
\\[-6pt]
\text{s. t.} \quad \forall j\in\{1,\ldots,m\},\quad
\left|u-\tilde{\alpha}_{j}\right| \leq \sqrt{r(j,\varepsilon)}.
\end{array}
\right.
$$
As a consequence,
$$ u_{j} = \left<e_{j},\hat{\alpha}^{CS}\right>
= sgn \left(\tilde{\alpha}_{j}\right)
\left(\left|\tilde{\alpha}_{j}\right|-\sqrt{r(j,\varepsilon)}\right)_{+} $$
that does not depend on $p$.
Note that $u_{j}$ is a thresholded estimator of the correlation between $Y$ and $f_{j}(X)$,
this is what suggested the name ``Correlation Selector''.
Let us put $U$ the column vector that contains the $u_{j}$ for $j\in\{1,\ldots,m\}$
and $M$ the matrix $(\left<e_{i},e_{j}\right>_{GN})_{i,j}$, then $\hat{\alpha}^{CS}$ is just any solution
of $\hat{\alpha}^{CS} M = U$.

\begin{rmk}
\label{rmkspar}
Note that the Correlation Selector has no reason to be sparse, however, the vector $\hat{\alpha}^{CS} M$
is sparse. An interpretation of this fact is given in the next subsection.
\end{rmk}

\subsection{Oracle inequality for the Correlation Selector}

\begin{thm}
\label{thm4}
We have:
$$ P\left[\left\|\hat{\alpha}^{CS}-\overline{\alpha}\right\|_{CS}^{2}
\leq \inf_{S\subset\{1,\ldots,m\}} \Bigg( \sum_{j\notin S}
\left<\overline{\alpha},e_{j}\right>_{GN}^{2} +4 \sum_{j\in S}
r(j,\varepsilon) \Bigg) \right] \geq 1-\varepsilon. $$ Moreover, if
we assume that there is a $D>0$ such that for any
$\alpha\in\mathcal{E}_{m}$, $\left\|\alpha\right\|_{GN} \geq D
\left\|\alpha\right\|$ where
$$ \mathcal{E}_{m} = \left\{\alpha\in\mathds{R}^{m},\quad \left<\overline{\alpha},e_{j}\right>_{GN} = 0
\Rightarrow \left<\alpha,e_{j}\right> = 0 \right\} $$
then we have:
$$ P\left[ \left\|\hat{\alpha}^{CS}-\overline{\alpha}\right\|_{GN}^{2}
\leq \frac{1}{D^{2}}\inf_{S\subset\{1,\ldots,m\}} \Bigg( \sum_{j\notin
S} \left<\overline{\alpha},e_{j}\right>_{GN}^{2} +4\sum_{j\in S}
r(j,\varepsilon) \Bigg) \right] \geq 1-\varepsilon. $$
\end{thm}

The proof can be found in Subsection \ref{proofthm4} page
\pageref{proofthm4}.

\begin{rmk}
Note that the result on
$\left\|\hat{\alpha}^{CS}-\overline{\alpha}\right\|_{CS}$ does not
require any assumption on the dictionary of functions. However, this
quantity does not have, in general, an interesting interpretation.
The result about the quantity of interest,
$\left\|\hat{\alpha}^{CS}-\overline{\alpha}\right\|_{GN}^{2}$,
requires that a part of the dictionary is almost orthogonal, this
condition is to be compared to the one in Theorem \ref{thTSYBA}.
\end{rmk}

\begin{rmk}
Note that if there is a $\overline{S}$ such that for any $j\notin \overline{S}$,
$\left<\overline{\alpha},e_{j}\right>_{GN} = 0$
and if $r(j,\varepsilon)=k\log(m/\varepsilon)/n$ then we have:
$$
P\left[\left\|\hat{\alpha}^{CS}-\overline{\alpha}\right\|_{CS}^{2}
\leq \frac{4k|\overline{S}|\log\frac{m}{\varepsilon}}{n} \right]
\geq 1-\varepsilon,
$$
and if moreover for any $\alpha\in\mathcal{E}_{m}$, $\left\|\alpha\right\|_{GN}
\geq D \left\|\alpha\right\|$
then
$$ P\left[ \left\|\hat{\alpha}^{CS}-\overline{\alpha}\right\|_{GN}^{2}
\leq \frac{4k|\overline{S}|\log\frac{m}{\varepsilon}}{D^{2}n}
\right] \geq 1-\varepsilon. $$ The condition that for a lot of $j$,
$\left<\overline{\alpha},e_{j}\right>_{GN} = 0$ means that most of
the functions in the dictionary are not correlated with $Y$. In
terms of sparsity, it means that the vector $\overline{\alpha} M$ is
sparse. So, intuitively, the Correlation Selector will perform well
when most of the functions in the dictionary have weak correlation
with $Y$, but we expect that altogether these functions can bring a
reasonable prediction for $Y$.
\end{rmk}

\section{Numerical simulations}

\label{simulations}

\subsection{Motivation}

We compare here LASSO, Iterative Feature Selection and Correlation
Selector on a toy example, introduced by Tibshirani \cite{Lasso1}.
We also compare their performances to the ordinary least square
(OLS) estimate as a benchmark. Note that we will not propose a very
fine choice for the $r(j,\varepsilon)$. The idea of these
simulations is not to identify a good choice for the penalization in
practice. The idea is to observe the similarity and differences
between different order in projections in our general algorithm,
using the same confidence regions.

\subsection{Description of the experiments}

The model defined by Tibshirani \cite{Lasso1} is the following. We have:
$$ \forall i\in\{1,\ldots,20\}, \quad Y_{i} = \left<\beta,X_{i}\right> + \varepsilon_{i} $$
with $X_{i}\in\mathcal{X}=\mathds{R}^{8}$, $\beta\in\mathds{R}^{8}$ and the $\varepsilon_{i}$
are i. i. d. from a gaussian distribution with mean $0$ and standard deviation $\sigma$.

The $X_{i}$'s are i. i. d. too, and each $X_{i}$ comes from a gaussian distribution with mean
$(0,\ldots,0)$ and with variance-covariance matrix:
$$\Sigma(\rho) = \bigl(\rho^{|i-j|}\bigr)_{
\tiny{
\begin{array}{c}
i\in\{1,\ldots,8\}
\\
j\in\{1,\ldots,8\}
\end{array}
}
}
$$
for $\rho\in[0,1[$.

We will use the three particular values for $\beta$ taken by Tibshirani \cite{Lasso1}:
\begin{align*}
\beta^{1} & = (3,1.5,0,0,2,0,0,0), \\
\beta^{2} & = (1.5,1.5,1.5,1.5,1.5,1.5,1.5,1.5), \\
\beta^{3} & = (5,0,0,0,0,0,0,0),
\end{align*}
corresponding to a ``sparse'' situation ($\beta^{1}$), a ``non-sparse'' situation ($\beta^{2}$)
and a ``very sparse'' situation ($\beta^{3}$).

We use two values for $\sigma$: 1 (the ``low noise case'') and 3 (the ``noisy case'').

Finally, we use two values for $\rho$: 0.1 (``weakly correlated variables'') and 0.5 (``highly correlated variables'').

We run each example (corresponding to a given value of $\beta$,
$\sigma$ and $\rho$) 250 times. We use the software R \cite{R} for
simulations. We implement Iterative Feature Selection as described
in subsection \ref{IFSdes} page \pageref{IFSdes}, and the
Correlation Selector, while using the standard OLS estimate and the
LASSO estimator given by the LARS package described in \cite{LARS}.
Note that we use the estimators defined in the {\bf deterministic
design case}, this means that we consider $\|.\|_{GN}=\|.\|_{n}$
(the empirical norm) as our criterion here. The choice:
$$ r(\varepsilon) = r(j,\varepsilon) = \frac{\sigma}{3}\sqrt{\frac{\log m}{n}} = \frac{\sigma}{3} \sqrt{\frac{\log 8}{20}} $$
was not motivated by theoretical considerations but seems to perform well in practice.

\subsection{Results and comments}

The results are reported in Table \ref{table1}.

\begin{table}[t!]
\caption{Results of Simulations.
For each possible combination of $\beta$, $\sigma$ and $\rho$,
we report in a column the mean empirical loss over the 250 simulations, the standard deviation
of this quantity over the simulations and finally the mean number of non-zero coefficients
in the estimate, for ordinary least square (OLS), LASSO, Iterative Feature
Selection (IFS) and Correlation Selector (C-SEL)}\label{table1}
\begin{tabular}{|p{2cm}|p{0.5cm}|p{0.5cm}||p{1.0cm}|p{1.0cm}|p{1.0cm}|p{1.0cm}|}
\hline
$\beta$ & $\sigma$ & $\rho$ & {OLS} & {LASSO} & {IFS} & {C-SEL} \\
\hline \hline
                     & $3$ & $0.5$ & \bf{3.67} & \bf{1.64} & \bf{1.56} & \bf{3.65} \\
$\beta^{1}$          &     &       & 1.84 & 1.25 & 1.20 & 1.96 \\
(sparse)             &     &       & 8 & 4.64 & 4.62 & 8 \\
\hline
                     & $1$ & $0.5$ & \bf{0.40} & \bf{0.29} & \bf{0.36} & \bf{0.44} \\
                     &     &       & 0.22 & 0.19 & 0.23 & 0.23 \\
                     &     &       & 8 & 5.42 & 5.70 & 8 \\
\hline
                     & $3$ & $0.1$ & \bf{3.75} & \bf{2.72} & \bf{2.85} & \bf{3.44} \\
                     &     &       & 1.86 & 1.50 & 1.58 & 1.72 \\
                     &     &       & 8 & 5.70 & 5.66 & 8 \\
\hline
                     & $1$ & $0.1$ & \bf{0.40} & \bf{0.30} & \bf{0.31} & \bf{0.43} \\
                     &     &       & 0.19 & 0.19 & 0.19 & 0.20 \\
                     &     &       & 8 & 5.92 & 5.96 & 8 \\
\hline \hline
                     & $3$ & $0.5$ & \bf{3.54} & \bf{3.36} & \bf{4.90} & \bf{3.98} \\
$\beta^{2}$          &     &       & 1.82 & 1.64 & 1.58 & 1.85 \\
(non sparse)         &     &       & 8 & 7.08 & 6.57 & 8 \\
\hline
                     & $1$ & $0.5$ & \bf{0.41} & \bf{0.54} & \bf{0.84} & \bf{0.47} \\
                     &     &       & 0.21 & 0.93 & 0.36 & 0.24 \\
                     &     &       & 8 & 7.94 & 7.89 & 8 \\
\hline
                     & $3$ & $0.1$ & \bf{3.78} & \bf{3.82} & \bf{4.50} & \bf{4.01} \\
                     &     &       & 1.78 & 1.51 & 1.59 & 1.86 \\
                     &     &       & 8 & 7.06 & 7.03 & 8 \\
\hline
                     & $1$ & $0.1$ & \bf{0.40} & \bf{0.42} & \bf{0.71} & \bf{0.48} \\
                     &     &       & 0.20 & 0.29 & 0.32 & 0.22 \\
                     &     &       & 8 & 7.98 & 7.98 & 8 \\
\hline \hline
                     & $3$ & $0.5$ & \bf{3.55} & \bf{1.65} & \bf{1.59} & \bf{3.42} \\
$\beta^{3}$          &     &       & 1.79 & 1.28 & 1.27 & 1.74 \\
(very sparse)        &     &       & 8 & 4.48 & 4.49 & 8 \\
\hline
                     & $1$ & $0.5$ & \bf{0.40} & \bf{0.18} & \bf{0.17} & \bf{0.46} \\
                     &     &       & 0.21 & 0.14 & 0.14 & 0.25 \\
                     &     &       & 8 & 4.46 & 4.48 & 8 \\
\hline
                     & $3$ & $0.1$ & \bf{3.46} & \bf{1.69} & \bf{1.62} & \bf{3.00} \\
                     &     &       & 1.74 & 1.29 & 1.18 & 1.45 \\
                     &     &       & 8 & 4.92 & 4.92 & 8 \\
\hline
                     & $1$ & $0.1$ & \bf{0.40} & \bf{0.20} & \bf{0.19} & \bf{0.44}\\
                     &     &       & 0.20 & 0.14 & 0.14 & 0.24 \\
                     &     &       & 8 & 4.98 & 4.91 & 8 \\
\hline
\end{tabular}
\end{table}

The following remarks can easily be made in view of the results:
\begin{itemize}
\item both methods based on projection on random confidence regions using the norm $\|.\|_{GN}=\|.\|_{n}$
clearly outperform the OLS in the
sparse cases, moreover they present the advantage of giving sparse estimates;
\item in the non-sparse case, the OLS performs generally better than the other methods, but LASSO is very close,
it is known that a better choice for the value $r(j,\varepsilon)$ would lead to a better result
(see Tibshirani \cite{Lasso1});
\item LASSO seems to be the best method on the whole set of experiments. In every case, it is never the worst method,
and always performs almost as well as the best method;
\item in the ``sparse case'' ($\beta^{1}$), note that IFS and LASSO are very close for the
small value of $\rho$. This is coherent with the previous theory, see
Theorem \ref{thm3} page \pageref{thm3};
\item IFS gives very bad results in the non-sparse case ($\beta^{2}$), but is the best method
in the sparse case ($\beta^{3}$). This last point tends to indicate that different situations
should lead to a different choice for the confidence regions we are to project on. However,
theoretical results leading on that choice are missing;
\item the Correlation Selector performs badly on the whole set of experiments. However,
note that the good performances for LASSO and IFS occurs for sparse values of $\beta$,
and the previous theory ensures good performances for C-SEL when $\beta M$ is sparse
where $M$ is the covariance matrix of the $X_{i}$. In other words, two experiments where
favorable to LASSO and IFS, but there was no experiment favorable to C-SEL.
\end{itemize}

In order to illustrate this last point, we build a new experiment favorable to C-SEL.
Note that we have
\begin{equation}
\label{exp1}
Y_{i} = \left<X_{i},\beta\right>  + \varepsilon_{i} = \left<X_{i}M^{-1},\beta M\right> + \varepsilon_{i}
\end{equation}
where $M$ is the correlation matrix of the $X_{i}$. Let us put $\tilde{X}_{i}= X_{i}M^{-1}$
and $\tilde{\beta} = \beta M$, we have the following linear model:
\begin{equation}
\label{exp2}
Y_{i} = \big<\tilde{X}_{i},\tilde{\beta}\big> + \varepsilon_{i} .
\end{equation}
The sparsity of $\beta$ gives advantage to the LASSO for estimating $\beta$ in Model
\ref{exp1}, it also gives an advantage to C-SEL for estimating $\tilde{\beta}$ in Model
\ref{exp2} (according to Remark \ref{rmkspar} page \pageref{rmkspar}).

We run again the experiments with $\beta=\beta^{3}$ and this time we try to estimate
$\tilde{\beta}$ instead of $\beta$ (so we act as if we had observed $\tilde{X}_{i}$
and not $X_{i}$).

Results are given in Table \ref{table2}.

\begin{table}[t!]
\caption{Results for the estimation of $\tilde{\beta}$.
As previously, for each possible combination of $\sigma$ and $\rho$,
we report in a column the mean empirical loss over the 250 simulations, the standard deviation
of this quantity over the simulations and finally the mean number of non-zero coefficients
in the estimate, this
for each estimate: OLS, LASSO, IFS and C-SEL}\label{table2}
\begin{tabular}{|p{2cm}|p{0.5cm}|p{0.5cm}||p{1.0cm}|p{1.0cm}|p{1.0cm}|p{1.0cm}|}
\hline
$\beta$ & $\sigma$ & $\rho$ & {OLS} & {LASSO} & {IFS} & {C-SEL} \\
\hline \hline
                     & $3$ & $0.5$ & \bf{3.64} & \bf{4.83} & \bf{5.12} & \bf{2.41} \\
$\beta^{1}$          &     &       & 1.99 & 2.53 & 2.64 & 1.92 \\
(sparse)             &     &       & 8 & 5.98 & 6.05 & 8 \\
\hline
                     & $1$ & $0.5$ & \bf{0.41} & \bf{1.09} & \bf{0.92} & \bf{0.26} \\
                     &     &       & 0.21 & 1.72 & 0.48 & 0.19 \\
                     &     &       & 8 & 7.11 & 7.40 & 8 \\
\hline
                     & $3$ & $0.1$ & \bf{3.65} & \bf{3.71} & \bf{3.72} & \bf{2.09} \\
                     &     &       & 1.71 & 1.96 & 1.99 & 1.40 \\
                     &     &       & 8 & 6.25 & 6.28 & 8 \\
\hline
                     & $1$ & $0.1$ & \bf{0.40} & \bf{0.47} & \bf{0.55} & \bf{0.23} \\
                     &     &       & 0.20 & 0.25 & 0.16 & 0.27 \\
                     &     &       & 8 & 7.35 & 7.38 & 8 \\
\hline
\end{tabular}
\vspace*{-6pt}
\end{table}

The correlation selector clearly outperforms the other methods in
this case.

\section{Conclusion}

\label{secconc}

\subsection{Comments on the results of the paper}

This paper provides a simple interpretation of well-known algorithms
of statistical learning theory in terms of orthogonal projections on
confidence regions. This very intuitive approach also provides tools
to prove oracle inequalities.

Simulations shows that methods based on confidence regions clearly
outperforms the OLS estimate in most examples. Actually, the
theoretical results and the experiments lead to the following
conclusion: in the case where we think that $\overline{\alpha}$ is
sparse, that means, if we assume that only a few functions in the
dictionary are relevant, we should use the LASSO or the Dantzig
Selector (we know that these estimators are almost equivalent since
\cite{Lasso3}); IFS can be seen as a good algorithmic approximation
of the LASSO in the orthogonal case. In the other cases, we should
think of another method of approximation (LARS, relaxed greedy
algorithm\ldots). When $\overline{\alpha}M$ is sparse, i. e. almost all
the functions in the dictionary are uncorrelated with $Y$, then we
the Correlation Selector seems to be a reasonable choice. This is, in
some way, the ``desperate case'', where for example for various reason
a practitioner thinks that he has the good set of variables to
explain $Y$, but he realizes that only a few of them are correlated
with $Y$ and that methods based on the selection of a small subset
of variables (LASSO, \ldots) leads to unsatisfying results.

\subsection{Extentions}

First, note that all the results given here in the deterministic
design case ($\|.\|_{GN}=\|.\|_{n}$) and in the random design case
($\|.\|_{GN}=\|.\|_{X}$) can be extended to another kind of
regression problem: the transductive case, introduced by Vapnik
\cite{Vapnik}. In this case, we assume that $m$ more pairs
$(X_{n+1},Y_{n+1})$,\ldots, $(X_{n+m},Y_{n+m})$ are drawn (i. i. d.
from $\mathds{P}$), and that $X_{n+1}$,\ldots, $X_{n+m}$ are given to
the statistician, whose task is now to predict the missing values
$Y_{n+1}$,\ldots, $Y_{n+m}$. Here, we can introduce the following
criterion
$$ \left\|\alpha-\alpha'\right\|_{trans}^{2} = \frac{1}{m}
\sum_{i=n+1}^{m}\Bigg[\sum_{j=1}^{m}\alpha_{j}f_{j}(X_{i})-\sum_{j=1}^{m}\alpha_{j}'f_{j}(X_{i}) \Bigg]^{2} .$$
In \cite{Alquier2007}, we argue that this case is of considerable interest in practice, and we show that
Assumption {\bf (CRA)} can be satisfied in this context. So, the reader can check that all the results in the
paper can be extended to the case $\|.\|_{GN}=\|.\|_{trans}$.

Also note that this approach can easily be extended into general
statistical problems with quadratic loss: in our paper
\cite{AlquierDens}, the Iterative Feature Selection method is
generalized to the density estimation with quadratic loss problem,
leading to a proposition of a LASSO-like program for density
estimation, that have also been proposed and studied by Bunea,
Tsybakov and Wegkamp \cite{SPADES} under the name SPADES.

\subsection{Future works}

Future works on this topic include a general study of the projection into the intersection
of the confidence regions
$$
\left\{
\begin{array}{l}
\arg\min_{\alpha=(\alpha_{1},\ldots,\alpha_{m})\in\mathds{R}^{m}} \delta(\alpha,0)
\\
\\[-6pt]
\text{s. t.} \quad \forall j\in\{1,\ldots,m\},\quad
\left|\left<\alpha,e_{j}\right>_{GN}-\tilde{\alpha}_{j}\right| \leq
\sqrt{r(j,\varepsilon)}
\end{array}
\right.
$$
for a generic distance $\delta(.,.)$.

A generalization to confidence regions defined by grouped variables, that would include the Group
LASSO studied by Bakin \cite{Bakin}, Yuan and Lin
\cite{GLasso} and Chesneau and Hebiri \cite{GVLasso} as a particular case is also feasible.

A more complete experimental study, including comparison of various choices for $\delta(.,.)$ and for $r(j,\varepsilon)$
based {\it on theoretical results or on heuristics} would be of great interest.

\section{Proofs}

\label{proofs}

\subsection[Proof of Proposition 3.1]{Proof of Proposition \ref{PropLASSO}}

\label{ProofLASSO}

\begin{proof}
Let us remember program \ref{ProgLASSO}:
\begin{equation}
\left\{
\begin{array}{l}
\max_{\alpha\in\mathds{R}^{m}} - \|\alpha\|_{GN}^{2}
\\
\\[-6pt]
\text{s. t.} \quad \forall j\in\{1,\ldots,m\},\quad
\left|\left<\alpha,e_{j}\right>_{GN}-\tilde{\alpha}_{j}\right| \leq
\sqrt{r(j,\varepsilon)}.
\end{array}
\right.
\end{equation}
Let us write the lagrangian of this program:
\begin{multline*}
\mathcal{L}(\alpha,\lambda,\mu) = - \sum_{i}\sum_{j}\alpha_{i}\alpha_{j} \left<e_{i},e_{j}\right>_{GN}
\\
     + \sum_{j}\lambda_{j} \Bigg[\sum_{i}\alpha_{i} \left<e_{i},e_{j}\right>_{GN}
              - \tilde{\alpha}_{j} - \sqrt{r(j,\varepsilon)}\Bigg]
\\
     + \sum_{j}\mu_{j} \Bigg[-\sum_{i}\alpha_{i} \left<e_{i},e_{j}\right>_{GN}
              + \tilde{\alpha}_{j} - \sqrt{r(j,\varepsilon)}\Bigg]
\end{multline*}
with, for any $j$, $\lambda_{j}\geq 0$, $\mu_{j} \geq 0$ and $\lambda_{j}\mu_{j}=0$.
Any solution $(\alpha^{*})$ of Program \ref{ProgLASSO} must satisfy, for any $j$,
$$
0 = \frac{\partial \mathcal{L}}{\partial \alpha_{j}}(\alpha^{*},\lambda,\mu)
= - 2 \sum_{i} \alpha_{i}^{*} \left<e_{i},e_{j}\right>_{GN}
     + \sum_{i} \left(\lambda_{i} - \mu_{i}\right) \left<e_{i},e_{j}\right>_{GN},
$$
so for any $j$,
\begin{equation}
\label{SOLUTION}
 \sum_{i} \left<\frac{1}{2}\left(\lambda_{i} - \mu_{i}\right) e_{i},e_{j}\right>_{GN}
         =     \left<\alpha^{*},e_{j}\right>_{GN}.
\end{equation}
Note that this also implies that:
\begin{align*}
\|\alpha^{*}\|_{X} ={}& \bigg<\sum_{i}\alpha^{*}_{i}e_{i},\sum_{j}\alpha^{*}_{j}e_{j}\bigg>_{GN}
                   = \sum_{i}\alpha^{*}_{i} \bigg<e_{i},\sum_{j}\alpha^{*}_{j}e_{j}\bigg>_{GN}
\\
={}& \sum_{i}\alpha^{*}_{i} \bigg<e_{i},\sum_{j}\frac{1}{2}(\lambda_{j} - \mu_{j})e_{j}\bigg>_{GN}
= \sum_{j}\frac{1}{2}(\lambda_{j} - \mu_{j}) \bigg<\sum_{i}\alpha^{*}_{i}e_{i},e_{j}\bigg>_{GN}
\\
={}& \sum_{j}\sum_{i} \frac{1}{2}
(\lambda_{j} - \mu_{j}) \frac{1}{2}(\lambda_{i} - \mu_{i})  \left<e_{i},e_{j}\right>_{GN}.
\end{align*}
Using these relations, the lagrangian may be written:
\begin{multline*}
\mathcal{L}(\alpha^{*},\lambda,\mu)
= - \sum_{i}\sum_{j}\frac{1}{2}(\lambda_{i}-\mu_{i})
\frac{1}{2}(\lambda_{j}-\mu_{j}) \left<e_{i},e_{j}\right>_{GN}
\\
 + \sum_{i}\sum_{j}\frac{1}{2}(\lambda_{i}-\mu_{i})
      (\lambda_{j}-\mu_{j}) \left<e_{i},e_{j}\right>_{GN}
\\
 - \sum_{j}(\lambda_{j}-\mu_{j})\tilde{\alpha}_{j}
        +  \sum_{j} (\lambda_{j}+\mu_{j})
   \sqrt{r(j,\varepsilon)}.
\end{multline*}
Note that the condition $\lambda_{j}\geq 0$, $\mu_{j} \geq 0$ and $\lambda_{j}\mu_{j}=0$
means that there is a $\gamma_{j}\in\mathds{R}$ such that $\gamma_{j}=2(\lambda_{j}-\mu_{j})$,
$|\gamma_{j}| = 2(\lambda_{j}+\mu_{j})$, and so $\mu_{j}=(\gamma_{j}/2)_{-}$
and $\lambda_{j}=(\gamma_{j}/2)_{+}$. Let also $\gamma$ denote the vector which $j$-th component
is exactly $\gamma_{j}$, we obtain:
$$
\mathcal{L}(\alpha^{*},\lambda,\mu)
=
\left\|\gamma \right\|_{GN}^{2} - 2\sum_{j} \gamma_{j} \tilde{\alpha}_{j}
+  2 \sum_{j} \left|\gamma_{j}\right|
   \sqrt{r(j,\varepsilon)}
$$
that is maximal with respect to the $\lambda_{j}$ and $\mu_{j}$, so with respect
to $\gamma$. So $\gamma$ is a solution of Program \ref{lassogeneral}.

Now, note that Equation \ref{SOLUTION} ensures that any solution $\alpha^{*}$ of Program
\ref{ProgLASSO} satisfies:
\begin{equation*}
\bigg<\sum_{i} \gamma_{i} e_{i},e_{j}\bigg>_{GN}
         =   \left<\alpha^{*},e_{j}\right>_{GN}.
\end{equation*}
We can easily see that $\alpha^{*}=\gamma$ is a possible solution.

In the case where $\|.\|_{GN}$ is the empirical norm $\|.\|_{n}$ we obtain:
\begin{align*}
\|\gamma\|_{GN}^{2}-2\sum_{j=1}^{m}\gamma_{j}\tilde{\alpha}_{j} & =
\frac{1}{n}\sum_{i=1}^{n}\Bigg[\sum_{j=1}^{m}\gamma_{j}f_{j}(X_{i})\Bigg]^{2}
-2\frac{1}{n}\sum_{i=1}^{n} Y_{i}
\Bigg[\sum_{j=1}^{m}\gamma_{j}f_{j}(X_{i})\Bigg]
\\
& = \frac{1}{n}\sum_{i=1}^{n}\Bigg[Y_{i}
-\sum_{j=1}^{m}\gamma_{j}f_{j}(X_{i})\Bigg]^{2} -
\frac{1}{n}\sum_{i=1}^{n}Y_{i}^{2}.
\end{align*}
\end{proof}

\subsection[Proof of Theorem 3.2]{Proof of Theorem \ref{thm3}}

\label{proofthm3}

\begin{proof}
In the case of orthogonality, we have $\|.\|_{GN}=\|.\|$ the euclidian norm. So $\hat{\alpha}^{L}$
satisfies, according to its definition:
$$
\left\{
\begin{array}{l}
\displaystyle\arg\min_{\alpha=(\alpha_{1},\ldots,\alpha_{m})\in\mathds{R}^{m}}
\sum_{j=1}^{m}\alpha_{j}^{2}
\\
\\[-9pt]
\text{s. t.} \quad \forall j\in\{1,\ldots,m\},\quad
\left|\alpha_{j}-\tilde{\alpha}_{j}\right|\leq \sqrt{r(j,\varepsilon)}
\end{array}
\right.
$$
while $\hat{\alpha}^{DS}$ satisfies:
$$
\left\{
\begin{array}{l}
\displaystyle\arg\min_{\alpha=(\alpha_{1},\ldots,\alpha_{m})\in\mathds{R}^{m}}
\sum_{j=1}^{m}|\alpha_{j}|
\\
\\[-9pt]
\text{s. t.} \quad \forall j\in\{1,\ldots,m\},\quad
\left|\alpha_{j}-\tilde{\alpha}_{j}\right|\leq \sqrt{r(j,\varepsilon)}.
\end{array}
\right.
$$
We can easily solve both problem by an individual optimization on each $\alpha_{j}$ and obtain
the same solution
$$ \alpha_{j}^{*}
=sgn\left(\tilde{\alpha}_{j}\right) \left(\left|\tilde{\alpha}_{j}\right|-\sqrt{r(j,\varepsilon)}\right)_{+}.$$
For $\hat{\alpha}^{IFS}$ just note that in the case of orthogonality, sequential projections on each
$\mathcal{CR}(j,\varepsilon)$ leads to the same result than the projection on their intersection, so
$\hat{\alpha}^{IFS}=\hat{\alpha}^{L}$.
Then, let us choose $S\subset\{1,\ldots,m\}$ and remark that
\begin{multline*}
\left\|\hat{\alpha}^{L}-\overline{\alpha}\right\|_{GN}^{2}
= \left\|\hat{\alpha}^{L}-\overline{\alpha}\right\|^{2}
= \sum_{j=1}^{m} \left<\hat{\alpha}^{L}-\overline{\alpha},e_{j}\right>^{2}
\\
= \sum_{j\in S} \left<\hat{\alpha}^{L}-\overline{\alpha},e_{j}\right>^{2}
 + \sum_{j\notin S} \left<\hat{\alpha}^{L}-\overline{\alpha},e_{j}\right>^{2}.
\end{multline*}
Now, with assumption CRA, with probability $1-\varepsilon$, for any $j$,
$\overline{\alpha}$ satisfies the same constraint than the LASSO estimator so
$$\left| \left< \overline{\alpha},e_{j}\right>-\tilde{\alpha}_{j}\right| \leq \sqrt{r(j,\varepsilon)}$$
and so
$$
\left|\left<\hat{\alpha}^{L}-\overline{\alpha},e_{j}\right>\right|
= \left|\alpha_{j}^{*} - \left<\overline{\alpha},e_{j}\right> \right|
\leq \left|\alpha_{j}^{*} - \tilde{\alpha}_{j} \right|
     + \left|\left<\overline{\alpha},e_{j}\right> - \tilde{\alpha}_{j} \right|
\leq 2 \sqrt{r(j,\varepsilon)}.
$$
Moreover, let us remark that $\alpha_{j}^{*}$ is the number with the smallest
absolute value satisfying this contraint, so
$$
\left|\alpha_{j}^{*} - \left<\overline{\alpha},e_{j}\right> \right|
\leq \max\left(\left|\alpha_{j}^{*}\right|,\left|\left<\overline{\alpha},e_{j}\right> \right|\right)
\leq
\left|\left<\overline{\alpha},e_{j}\right> \right|.
$$
So we can conclude
$$
\left\|\hat{\alpha}^{L}-\overline{\alpha}\right\|_{GN}^{2}
\leq
\sum_{j\in S} 4 r(j,\varepsilon) + \sum_{j\notin S} \left<\overline{\alpha},e_{j}\right>^{2}
=
4 \sum_{j\in S} r(j,\varepsilon) + \left\|\overline{\alpha}-\overline{\alpha}_{S}\right\|^{2}.
$$
\end{proof}

\subsection[Proof of Theorem 3.4]{Proof of Theorem \ref{thm4}}

\label{proofthm4}

\begin{proof}
Note that, for any $S$:
\begin{multline*} \left\|\hat{\alpha}^{CS}-\overline{\alpha}\right\|_{CS}^{2} = \sum_{j=1}^{m}
\left<\hat{\alpha}_{CS}-\overline{\alpha},e_{j}\right>_{GN}^{2}
\\
= \sum_{j\in S}
\left<\hat{\alpha}^{CS}-\overline{\alpha},e_{j}\right>_{GN}^{2}
+ \sum_{j\notin S}
\left<\hat{\alpha}^{CS}-\overline{\alpha},e_{j}\right>_{GN}^{2}.
\end{multline*}
By the constraint satisfied by $\hat{\alpha}^{CS}$ we have:
$$ \left<\hat{\alpha}^{CS}-\overline{\alpha},e_{j}\right>_{GN}^{2} \leq
4 r(j,\varepsilon) .$$
Moreover, we must remember that $u_{j} = \left<\hat{\alpha}^{CS},e_{j}\right>_{GN}$ satisfies
the program
$$\left\{
\begin{array}{l}
\arg\min_{u} |u|
\\
\\[-6pt]
\text{s. t.} \quad
\forall j\in\{1,\ldots,m\},\quad
\left|u-\tilde{\alpha}_{j}\right| \leq
\sqrt{r(j,\varepsilon)},
\end{array}
\right.
$$
that is also satisfied by $\left<\overline{\alpha},e_{j}\right>_{GN}$,
so $|u_{j}| \leq |\left<\overline{\alpha},e_{j}\right>|$
and so
$$ \left|u_{j}-\left<\overline{\alpha},e_{j}\right>\right|
\leq \max \left(|u_{j}|,|\left<\overline{\alpha},e_{j}\right>|\right)
  = |\left<\overline{\alpha},e_{j}\right>| $$
and so we have the relation:
$$ \left<\hat{\alpha}^{CS}-\overline{\alpha},e_{j}\right>_{GN}^{2} \leq
  \left<\overline{\alpha},e_{j}\right>_{GN}^{2} . $$

So we obtain:
$$
 \left\|\hat{\alpha}^{CS}-\overline{\alpha}\right\|_{CS}^{2}
\leq
\sum_{j\in S} 4 r(j,\varepsilon)
+ \sum_{j\notin S}
\left<\overline{\alpha},e_{j}\right>_{GN}^{2}
$$
This proves the first inequality of the theorem. For the second one, we just have
to prove that $
(\hat{\alpha}^{CS}-\overline{\alpha})M \in\mathcal{E}_{m}$. But
this is trivial because of the relation:
$$ \left<(\hat{\alpha}^{CS}-\overline{\alpha}) M,e_{j}\right>^{2}
  = \left<\hat{\alpha}^{CS}-\overline{\alpha},e_{j}\right>_{GN}^{2} \leq
  \left<\overline{\alpha},e_{j}\right>_{GN}^{2} .$$
\end{proof}

\bibliographystyle{acm}

\end{document}